\documentclass{amsart}
\usepackage{amsmath, amssymb, amsfonts}
\usepackage{fullpage}
\usepackage{color}

\theoremstyle{plain}
\newtheorem{theorem}{Theorem}[section]

\numberwithin{equation}{section}

\theoremstyle{definition}
\newtheorem{definition}[theorem]{Definition}

\newcommand{\G}{\Gamma}
\newcommand{\Hyp}{\mathbb{H}}
\newcommand{\cA}{\mathcal{A}}
\newcommand{\cP}{\mathcal{P}}
\newcommand{\cL}{\mathcal{L}}
\newcommand{\F}{\mathbb{F}}
\newcommand{\Q}{\mathbb{Q}}

\newcommand{\X}{\mathbb{X}}

\DeclareMathOperator\Aut{Aut}
\DeclareMathOperator\Lk{Lk}
\DeclareMathOperator\Stab{Stab}
\DeclareMathOperator\CAT{CAT}

\DeclareMathOperator\SL{SL}

\def\polhk#1{\setbox0=\hbox{#1}{\ooalign{\hidewidth
    \lower1.0ex\hbox{$\,\lhook$}\hidewidth\crcr\unhbox0}}}
\newcommand{\Swiatkowski}{\'Swi{\polhk{a}}tkowski}

\title{Lattices in hyperbolic buildings}
\author{Anne Thomas}
\address{School of Mathematics and Statistics F07, University of Sydney NSW 2006, Australia}
\thanks{The author is supported in part by ARC Grant No. DP110100440.}

\begin{document}

\maketitle

\section*{Introduction}

This survey is intended as a brief introduction to the theory of hyperbolic buildings and their lattices.  Hyperbolic buildings are negatively curved geometric objects which also have a rich algebraic and combinatorial structure, and the study of these buildings and the lattices in their automorphism groups involves a fascinating mixture of techniques from many different areas of mathematics.

Roughly speaking, a hyperbolic building is obtained by gluing together many hyperbolic spaces which are tiled by polyhedra.  For the precise definition, together with background on general buildings and known constructions of hyperbolic buildings, see Section \ref{s:buildings} below.

Given a hyperbolic building $\Delta$, we write $G = \Aut(\Delta)$ for the group of automorphisms, or cellular isometries, of $\Delta$.  When the building $\Delta$ is locally finite, the group $G$ equipped with the compact-open topology is naturally a locally compact topological group, and so has a Haar measure $\mu$.  In this topology on $G$, a subgroup $\G < G$ is discrete if and only if it acts on $\Delta$ with finite cell stabilisers.  A \emph{lattice} in $G$ is a discrete subgroup $\G < G$ such that $\mu(\G \backslash G) < \infty$, and a lattice $\G$ is \emph{cocompact} (or \emph{uniform}) if $\G \backslash G$ is compact.  The Haar measure $\mu$ on $G$ may be normalised so that the covolume $\mu(\G \backslash G)$ of a lattice $\G < G$ is given by the formula
\begin{equation}\label{e:lattice}  \mu(\G \backslash G) = \sum \frac{1}{|\Stab_\G(v)|}  \end{equation}
where the sum is taken over a set of representatives for the orbits of the vertices of $\Delta$ under the action of $\Gamma$.  A discrete subgroup $\G < G$ is then a lattice if and only if this sum converges, and is a cocompact lattice if and only if $\Gamma \backslash G$ is compact.  

For brevity, we will refer to lattices in the automorphism groups of hyperbolic buildings as lattices in hyperbolic buildings.  In Section \ref{s:lattices} we describe known constructions of lattices in hyperbolic buildings, then discuss many different questions concerning such lattices.  Much of the study of lattices in hyperbolic buildings is motivated by the well-developed theory of lattices in semisimple Lie groups.  

We recommend that the reader consult also the survey Farb--Hruska--Thomas~\cite{MR2807842}, which discusses more general polyhedral complexes and their automorphism groups and lattices.  In order to avoid repetition, we have concentrated here on questions concerning lattices which are particularly pertinent to hyperbolic buildings, and/or where there has been progress since \cite{MR2807842} was written.  We also provide greater detail than \cite{MR2807842} on hyperbolic Coxeter groups and constructions of hyperbolic buildings.  

\subsection*{Acknowledgements}  The author thanks the organisers of ``Geometry, Topology and Dynamics in Negative Curvature" for the opportunity to attend such a well-run and interesting conference, the London Mathematical Society for travel support and an anonymous referee for helpful comments.  The author is indebted to her coauthors on \cite{MR2807842} and to many of the researchers cited below for numerous rewarding discussions.

\section{Buildings and hyperbolic buildings}\label{s:buildings}

In this section we recall definitions and results concerning both general and hyperbolic buildings.  We begin with a summary of the relevant theory of Coxeter groups and Coxeter polytopes in Section \ref{ss:coxeter}, and some background on polyhedral complexes in Section \ref{ss:polyhedral}.  Buildings and examples of spherical and Euclidean buildings are discussed in Section \ref{ss:buildings} before we focus on hyperbolic buildings in Section \ref{ss:hyperbolic buildings}.  Some references for the theory of buildings are Abramenko--Brown~\cite{MR2439729}, Brown~\cite{MR1644630} and Ronan~\cite{MR2560094}.  

\subsection{Coxeter groups and Coxeter polytopes}\label{ss:coxeter}

We mostly follow the reference Davis~\cite{MR2360474}, particularly Chapter 6, and concentrate on the hyperbolic case.

Recall that a \emph{Coxeter group} is a group $W$ with finite generating set $S$ and presentation of the
form \[ W = \langle s\in S \mid (s t )^{m_{st}} = 1\rangle \] where $s, t \in S$, $m_{ss} = 1$ for all $s \in S$ and if
$s \neq t$ then $m_{st}$ is an integer $\geq 2$ or $m_{st} = \infty$, meaning that the product $st$ has infinite order.  The pair $(W,S)$ is called a \emph{Coxeter system}.  A Coxeter system $(W,S)$ is \emph{right-angled} if for each $s, t \in S$ with $s \neq t$, $m_{st} \in \{ 2, \infty\}$.  Note that $m_{st} = 2$ if and only if $st = ts$.

Let $\X^n$ be the $n$--dimensional sphere, $n$--dimensional Euclidean space or $n$--dimensional (real) hyperbolic space.  Many important examples of Coxeter groups arise as discrete reflection groups acting on $\X^n$, as follows.  Let $P$ be a convex polyhedron in $\X^n$ with all dihedral angles integer submultiples of $\pi$.  Such a $P$ is called a \emph{Coxeter polytope}.  Let $W = W(P)$ be the group generated by the set $S = S(P)$ of reflections in the codimension one faces of $P$.  Then $(W,S)$ is a Coxeter system and $W$ is a discrete subgroup of the isometry group of $\X^n$ (see~\cite[Theorem 6.4.3]{MR2360474}).  Moreover, the action of $W$ tessellates $\X^n$ by copies of $P$.  For example, let $P$ be a right-angled hyperbolic $p$--gon, $p \geq 5$.  Then the corresponding Coxeter system is right-angled with $p$ generators, one for each side of $P$, so that if $s$ and $t$ are reflections in distinct sides then $m_{st} = 2$ when these sides are adjacent, and otherwise $m_{st} = \infty$.

Suppose that $\X^n$ is the sphere or Euclidean space.  Then Coxeter polytopes $P \subset \X^n$ exist and have been classified in every dimension, and the corresponding Coxeter systems $(W,S)$ are the spherical or affine Coxeter systems, respectively (see~\cite[Table 6.1]{MR2360474} for the classification).  

If $\X^n$ is $n$--dimensional hyperbolic space $\Hyp^n$, then there is no complete classification of Coxeter polytopes.  Vinberg's Theorem~\cite{MR774946} establishes that compact hyperbolic Coxeter polytopes can exist only in dimension $n \leq 29$, although at the time of writing the highest dimension in which an example is known is $n = 8$ (due to Bugaenko~\cite{MR1155663}).  Finite volume hyperbolic Coxeter polytopes have also been investigated, with for example Prokhorov~\cite{MR842588} proving these can exist only in dimension $n \leq 995$.
For $n \leq 6$ (respectively, $n \leq 19$), there are infinitely many essentially distinct  compact (respectively, finite volume) hyperbolic Coxeter polytopes (Allcock~\cite{MR2240904}).   In dimension $3$, Andreev's Theorem~\cite{MR0259734} classifies compact hyperbolic Coxeter polytopes, but in dimensions $n \geq 4$ only special cases have been considered, and there seems little hope of a complete list.   Hyperbolic Coxeter polytopes which are simplices exist in dimensions $n \leq 4$ only, and their classification is given in \cite[Table 6.2]{MR2360474}.  Right-angled compact hyperbolic polytopes also exist in dimensions $n \leq 4$ only, and there are infinitely many examples in each dimension $n \leq 4$ (see \cite{MR1254933}).  A right-angled example in dimension $3$ is the dodecahedron, which tessellates $\Hyp^3$ as depicted on the cover of Thurston's book~\cite{MR1435975}, and a right-angled example in dimension $4$ is the $120$--cell, which has 120 dodecahedral faces.   For other special cases of compact hyperbolic Coxeter polytopes, see for example the work of Esselmann~\cite{MR1396674}, Felikson--Tumarkin~\cite{MR2378860,MR2513640,MR2549446}, Kaplinskaja~\cite{MR0360858} and Tumarkin~\cite{MR2350459,MR2086616,MR2193442}, and for an overview of results in the finite volume case, see the introduction to \cite{MR2240904}.

Let $W = W(P)$ be the Coxeter group generated by reflections in the faces of a hyperbolic Coxeter polytope $P$.  If $P$ is compact then $W$ is a word-hyperbolic group, that is, a group which is hyperbolic in the sense of Gromov.  (For background on word-hyperbolic groups, see \cite{MR1744486}.  Necessary and sufficient conditions for word-hyperbolicity of Coxeter groups were established by Moussong \cite[Corollary 12.6.3]{MR2360474}.)  On the other hand, some authors such as Humphreys~\cite{MR1066460} reserve ``hyperbolic Coxeter group" for the case that $P$ is a compact simplex.  In this survey, when necessary we will refer to the discrete reflection group $W(P)$, where $P$ is a (compact) hyperbolic Coxeter polytope, as a \emph{(cocompact) geometric hyperbolic Coxeter group}.  

\subsection{Polyhedral complexes and links}\label{ss:polyhedral}

Polyhedral complexes are generalisations of (geometric realisations of) simplicial complexes.  Roughly speaking, they are obtained by gluing together polyhedra from the constant curvature space $\X^n$ (the sphere, Euclidean space or hyperbolic space), using isometries along faces.  For the formal definition of a polyhedral complex, see for example \cite[Section 2.1]{MR2807842}.  We sometimes refer to $2$--dimensional polyhedral complexes as \emph{polygonal complexes}.

The tessellation of $\X^n$ by copies of a Coxeter polytope $P$ is a simple example of a polyhedral complex.  A metric tree is a $1$--dimensional Euclidean polyhedral complex, and a product of two such trees is a $2$--dimensional Euclidean polygonal complex.  

Let $x$ be a vertex of an $n$--dimensional polyhedral complex $X$. The \emph{link} of $x$, denoted $\Lk(x,X)$, is the spherical $(n-1)$--dimensional  polyhedral complex obtained by intersecting $X$ with an $n$--sphere of sufficiently small radius centred at $x$.  For example, if $X$ has dimension $2$, then $\Lk(x,X)$ may be identified with the graph whose vertices correspond to endpoints of edges of X that are incident to $x$, and whose edges correspond to corners of faces of $X$ incident to $x$.  By rescaling so that for each $x$ the $n$--sphere around $x$ has radius $1$, we induce a canonical metric on each link.  

The importance of links is that they provide a local condition for nonpositive or negative curvature of polyhedral complexes, using the following result which combines several theorems of Gromov.  For these theorems as well as background on the nonpositive curvature condition $\CAT(0)$, and the negative curvature condition $\CAT(-1)$, see \cite{MR1744486}.

\begin{theorem}[Gromov]\label{t:CAT1} Let $X$ be a contractible polyhedral complex of piecewise constant curvature~$\kappa$.  If $X$ has finitely many isometry types of cells, then $X$ is $\CAT(\kappa)$ if and only if for all vertices $x$ of $X$, the link $\Lk(x,X)$ is a $\CAT(1)$ space.  In particular, if $X$ is a contractible Euclidean (respectively, hyperbolic) polygonal complex with finitely many isometry types of cells, then $X$ is $\CAT(0)$ (respectively, $\CAT(-1)$) if and only if every embedded loop in the graph $\Lk(x,X)$ has length at least $2\pi$.\end{theorem}

An important special case of a polyhedral complex is a \emph{$(k,L)$--complex}, which is a polygonal complex in which each face is a regular $k$--gon, for $k \geq 3$ an integer, and the link at each vertex is a fixed finite graph $L$.  So long as $k$ and $L$ satisfy a simple condition, Ballmann--Brin~\cite{MR1279883} showed that a contractible $\CAT(0)$ $(k,L)$--complex may be constructed by a ``free" inductive process of adding $k$--gons to the previous stage.  Another construction of $(k,L)$--complexes for $k$ even is the special case of the Davis--Moussong complex described in \cite[Section 3.6]{MR2807842}.

\subsection{Buildings}\label{ss:buildings}

We will adopt the following ``geometric" definition of a building, which is the most appropriate for the hyperbolic case.  Other definitions, using simplicial complexes or chamber systems, may be found in, for example,~\cite{MR2439729} or~\cite{MR2560094}.

\begin{definition}\label{d:building} Let $\X^n$ be respectively the $n$--dimensional sphere, $n$--dimensional Euclidean space or $n$--dimensional hyperbolic space.  Let $P$ be a Coxeter polytope in $\X^n$ and let $(W,S)$ be the corresponding Coxeter system.   A respectively \emph{spherical, Euclidean or hyperbolic building of type $(W,S)$} is a polyhedral complex $\Delta$ equipped with a maximal family of subcomplexes, called \emph{apartments}, so that each apartment is isometric to the tessellation of $\X^n$ by copies of $P$, called \emph{chambers},
and so that: \begin{enumerate} \item any two chambers of $\Delta$ are contained in a common apartment; and \item for any
two apartments $\cA$ and $\cA'$, there exists an isometry $\phi:\cA \to \cA'$ which fixes $\cA \cap \cA'$.
\end{enumerate}  \end{definition}

The tessellation of a single copy of $\X^n$ by images of $P$ satisfies Definition \ref{d:building} and is sometimes called a \emph{thin} building.  We will mainly be interested in \emph{thick} buildings, those where there is ``branching", that is, where each codimension one face of each chamber is contained in at least three distinct chambers.  Thick buildings may be thought of as obtained by gluing together many copies of the same tessellation of $\X^n$.  A building is \emph{right-angled} if it is of type $(W,S)$ a right-angled Coxeter system.

We now discuss some important examples of spherical and Euclidean buildings.  Hyperbolic buildings will be considered in Section \ref{ss:hyperbolic buildings} below.

\subsubsection{Examples of spherical buildings}\label{sss:spherical}

A first example of a spherical building is the complete bipartite graph $K_{q,q}$, which is thick so long as $q \geq 3$.  The chambers are the edges of this  graph, metrised as quarter-circles, and the apartments are the embedded loops of length $2\pi$.  This is a (right-angled) building of type $(W,S)$ where 
\[  W = \langle s, t \mid s^2 = t^2 = (st)^2 = 1 \rangle \]
is the dihedral group of order $4$, acting on the circle.  

An example of a spherical building $\Delta$ which is not right-angled is the flag complex of the projective plane over the finite field $\F_q$ of order $q$, which may be constructed as follows.  Let $V$ be the vector space $\F_q \times \F_q \times \F_q$ over $\F_q$, let $\cP$ be the collection of one-dimensional subspaces of $V$ (the \emph{points} of the projective plane) and let $\cL$ be the collection of two-dimensional subspaces of $V$ (the \emph{lines}).  A point $p \in \cP$ is defined to be \emph{incident} to a line $l \in \cL$ if $p \subset l$.  The building $\Delta$ is then the bipartite graph with vertex set $\cP \sqcup \cL$ and edges corresponding to incidence.  See for instance \cite[Chapter 1, Example 3]{MR2560094} for the verification that $\Delta$ is indeed a building, of type 
\[ W =   \langle s, t \mid s^2 = t^2 = (st)^3 = 1 \rangle \]
the symmetric group on three letters, which is isomorphic to the dihedral group of order $6$, acting on the circle.  In particular, the apartments of $\Delta$ are embedded cycles of $6$ edges, and correspond to bases of $V$.  The \emph{standard apartment} is that corresponding to the standard basis for $V$.  

The structure of this spherical building is strongly connected with the structure of the group $G = \SL_3(\F_q)$, in a way that generalises to many other pairings of buildings with groups.  The action of $G$ on $V$ induces a natural action on the building $\Delta$.  The stabilisers of edges of $\Delta$ are the cosets of the upper-triangular subgroup $B < G$, the (standard) \emph{Borel subgroup}, and the stabilisers of vertices are the \emph{parabolic subgroups} of $G$.  Thus  the chambers of $\Delta$ can be identified with the cosets $G/B$, and the vertices of $\Delta$ can be identified with the disjoint union $G/P_1 \sqcup G/P_2$, where $P_1$ is the parabolic subgroup fixing the span of $(1,0,0)$ and $P_2$ that fixing the span of $(1,0,0)$ and $(0,1,0)$.  The pointwise stabiliser of the standard apartment is the diagonal subgroup $T < G$, called the \emph{torus}, and the setwise stabiliser of the standard apartment is the group of monomial matrices $N < G$, that is, matrices with exactly one non-zero entry in each row and each column.  The normaliser of the torus $T$ is the group $N$, with quotient $N/T \cong W$, and the group $W$ is called the \emph{Weyl group}.  

The discussion in the previous paragraph could be summarised by saying that the group $G = \SL_3(\F_q)$ has a \emph{$(B,N)$--pair}, also known as a \emph{Tits system}; there is then a building associated to $G$, of type its Weyl group, which is constructed using cosets of important subgroups as indicated.   (For the rather technical definition of a $(B,N)$--pair, see for example \cite[Chapter 5]{MR2560094}.)  Other spherical buildings are associated to other finite groups of Lie type, using their respective $(B,N)$--pairs.  

All one-dimensional spherical buildings are \emph{generalised $m$--gons}, meaning that they are graphs with diameter $m$ edges and shortest embedded circuit containing $2m$ edges.   A finite, thick generalised $m$--gon exists for $m \in \{ 2,3,4,6,8\}$ only (Feit--Higman, see~\cite[Theorem 3.4]{MR2560094}).  Generalised $2$--gons are complete bipartite graphs.  Generalised $3$--gons are flag complexes of projective planes, and so there is no classification known.  There are many examples of finite generalised $4$--gons, but only one or two known examples of finite generalised $6$-- or $8$--gons.  For more on generalised $m$--gons, see for example \cite[Chapter 5]{MR1829620}.

\subsubsection{Examples of Euclidean buildings}\label{sss:euclidean}

 A first example of a Euclidean building is the tree $T_q$ of valence $q$, metrised so that each edge has length say $1$.  The chambers are the edges of the tree, and the apartments are the bi-infinite geodesics.  This is a (right-angled) building of type $(W,S)$ where
\[ W = \langle s, t \mid s^2 = t^2 = 1 \rangle \]
is the infinite dihedral group acting on the real line, with the generating reflections $s$ and $t$ fixing points distance $1$ apart.  The product of trees $T_q \times T_q$ is a $2$--dimensional right-angled Euclidean building with apartments the tessellation of the Euclidean plane by unit squares, and associated Coxeter system the direct product of two infinite dihedral groups.  

Many Euclidean buildings (and all irreducible Euclidean buildings of dimension $\geq 3$) are of ``algebraic" origin, as in the following example.  Let $K$ be a nonarchimedean local field, such the $p$--adics $\Q_p$ or the field of formal Laurent series $\F_q((t))$.  The group $G = \SL_n(K)$ has a $(B,N)$--pair with associated Weyl group $W$ an affine Coxeter group, so that the action of $W$ tessellates $(n-1)$--dimensional Euclidean space.  Thus the group $G$ has an associated Euclidean building of dimension $(n-1)$.  For instance the building for $\SL_3(\Q_p)$ has apartments the tessellation of the Euclidean plane by equilateral triangles.  For further details and references concerning Euclidean buildings, which are also known as \emph{affine buildings}, see \cite[Section 3.1]{MR2807842}.

\subsubsection{Links of buildings}\label{sss:links buildings}

Let $x$ be a vertex of an $n$--dimensional building.  Then it is easy to verify that the link of $x$ is a spherical building of dimension $(n-1)$, with the induced apartment and chamber structure.  For example, the link of each vertex in $T_q \times T_q$ is the complete bipartite graph $K_{q,q}$, and the link of each vertex in the building for $\SL_3(\F_q((t)))$ is the spherical building $\Delta$ for $\SL_3(\F_q)$ from Section \ref{sss:spherical} above.  

With the natural piecewise spherical structure, a spherical building is a $\CAT(1)$ space.  This was shown by Davis \cite{MR1709955}, generalising a result of Gromov~\cite{MR919829} for right-angled spherical Coxeter systems and of Moussong for all spherical Coxeter systems (see \cite[Theorem 12.3.3]{MR2360474}).  By Theorem~\ref{t:CAT1} above, it follows that Euclidean (respectively, hyperbolic) buildings are $\CAT(0)$ spaces (respectively, $\CAT(-1)$ spaces).  The result that irreducible Euclidean buildings are $\CAT(0)$ was already well-known, and a proof can be found in \cite[Section 11.2]{MR2439729}.

\subsection{Hyperbolic buildings}\label{ss:hyperbolic buildings}

We first discuss examples and constructions of hyperbolic buildings.  Many of these constructions of hyperbolic buildings also yield lattices, as discussed in Section \ref{s:lattices} below.  We then briefly discuss the classification of hyperbolic buildings.

A hyperbolic building $\Delta$ of dimension $2$ is sometimes called a \emph{Fuchsian} building, since if $\Delta$ has type $(W,S)$ then $W$ may be regarded as a Fuchsian group.  By the restrictions on the dimension of hyperbolic Coxeter polytopes discussed in Section \ref{ss:coxeter} above, hyperbolic buildings with compact (respectively, finite volume) chambers can exist only in dimension $n \leq 29$ (respectively, $n \leq 995$), and right-angled hyperbolic buildings with compact chambers exist only in dimensions $2$, $3$ and $4$.

A first example of a hyperbolic building is \emph{Bourdon's building} $I_{p,q}$, defined and studied in~\cite{MR1445387}.  Let $P$ be a regular right-angled hyperbolic $p$--gon, with $p \geq 5$.  Then $P$ is a Coxeter polytope, and the building $I_{p,q}$ has type the associated right-angled Coxeter system.  Thus $I_{p,q}$ is a right-angled Fuchsian building with apartments hyperbolic planes tessellated by copies of $P$.  The link of each vertex of $I_{p,q}$ is the complete bipartite graph $K_{q,q}$, with $q \geq 2$, and so by Theorem \ref{t:CAT1} above, Bourdon's building is CAT$(-1)$.  Each edge of $I_{p,q}$ is contained in $q$ chambers, thus $I_{p,q}$ is thick for $q \geq 3$.   Bourdon's building can be thought of as a hyperbolic version of the product of trees $T_q \times T_q$.  However it  is not globally a product space.  
It is a $(k,L)$--complex with $k = p$ and $L = K_{q,q}$, and as descibed in~\cite{MR1445387} may be constructed using the Ballmann--Brin inductive process (see Section \ref{ss:polyhedral} above). 

A slightly more general example is the right-angled Fuchsian building $I_{p,\mathbf{q}}$ where $\mathbf{q} = (q_s)_{s \in S}$ is a $p$--tuple of integers $q_s \geq 2$, indexed by the generators of the Coxeter system $(W,S)$ associated to the regular right-angled $p$--gon $P$.  The edges of $I_{p,\mathbf{q}}$ are assigned types $s \in S$, and the vertices then inherit types $\{s,t\} \subset S$ with $m_{st} = 2$, so that the natural action of $W$ on each apartment is type-preserving.  The parameters $q_s$ record that each edge of type $s$ is contained in $q_s$ chambers, which makes the building \emph{regular}.  Each vertex of type $\{s,t\}$ has link the complete bipartite graph $K_{q_s, q_t}$. 

Any regular right-angled building may be constructed using complexes of groups.  (Note that not all right-angled buildings are hyperbolic, since for example a product of trees is not a hyperbolic building.)  For the general theory of complexes of groups see \cite{MR1744486}, and for a summary of this theory in the context of polyhedral complexes see \cite[Section 2.3]{MR2807842}.   The construction for $I_{p,\mathbf{q}}$ that we now sketch appears in Bourdon~\cite[Example 1.5.1(a)]{MR1756974}, was known earlier to Davis and Meier and is equivalent to a special case of constructions given in~\cite[Section 5]{MR1756974} and~\cite[Example 18.1.10]{MR2360474}.  Let the $p$--gon $P$ and parameters $\mathbf{q}$ be as in the previous paragraphs.  For each $s \in S$ let $G_s$ be a group of order $q_s$.  A complex of groups over $P$ with universal cover the building $I_{p,\mathbf{q}}$ is obtained by assigning groups to the face, edges and vertices of $P$ as follows.  The face group is trivial, the group on the edge of type $s$ is $G_s$ and the group on the vertex of type $\{s,t\}$ with $m_{st} = 2$ is the direct product $G_s \times G_t$.  Very roughly speaking, the universal cover is obtained by ``unfolding" this orbifold-like data.  

An example of a Fuchsian building which is not right-angled and is constructed using complexes of groups is as follows.  Let $P$ be a regular hyperbolic $k$--gon, with $k \geq 6$ even and all dihedral angles $\frac{\pi}{3}$.  The universal cover of the following complex of groups over $P$, from Gaboriau--Paulin \cite[Section 3.1]{MR1877215}, is a Fuchsian building with apartments hyperbolic planes tessellated by copies of $P$, and the link of every vertex the spherical building $\Delta$ described in Section \ref{sss:spherical} above.  Let $G = \SL_3(\F_q)$ and let $B$, $P_1$ and $P_2$ be the Borel and parabolic subgroups of $G$ as in Section \ref{sss:spherical} above.  The face group is $B$, the edge groups alternate between $P_1$ and $P_2$ and all vertex groups are $G$.  In \cite[Section 3.4]{MR1877215} Gaboriau--Paulin also use complexes of groups to construct some Fuchsian buildings with non-compact chambers.  In~\cite[Example 1.5.3]{MR1445387} Bourdon uses complexes of groups to construct Fuchsian buildings with right-angled triangles as chambers.

Some additional constructions of Fuchsian buildings are as follows.  Suppose $L$ is a one--dimensional spherical building and $k \geq 3$ is even.  Then a Davis--Moussong $(k,L)$--complex (see Section \ref{ss:polyhedral} above) may be metrised as a hyperbolic building with all links $L$ and all chambers $k$--gons.  Vdovina~\cite{MR1938699} constructed various Fuchsian buildings with even-sided chambers as universal covers of finite polygonal complexes whose links are one-dimensional spherical buildings, with not necessarily the same link at each vertex, while Kangaslampi--Vdovina~\cite{MR2665778} using similar techniques constructed Fuchsian buildings with chambers $n$--gons, $n \geq 3$, and links generalised $4$--gons.  Bourdon~\cite[Example 1.5.2]{MR1445387} obtained certain Fuchsian buildings by ``hyperbolising" affine buildings.

Recall from Sections~\ref{sss:spherical} and~\ref{sss:euclidean} above that some spherical and Euclidean buildings are obtained as buildings for groups which have $(B,N)$--pairs.  In similar fashion, some hyperbolic buildings arise as buildings for Kac--Moody groups.  A Kac--Moody group $\Lambda$ over a finite field $\F_q$ may be thought of as an infinite-dimensional analogue of an algebraic group over a nonarchimedean local field.  The group $\Lambda$ has twin $(B,N)$--pairs, which yield isomorphic twin buildings $\Delta_+$ and $\Delta_-$, with the group $\Lambda$ acting diagonally on the product $\Delta_+ \times \Delta_-$.  When the Weyl group $W$ of $\Lambda$ is a (cocompact) geometric hyperbolic Coxeter group, then the associated buildings $\Delta_\pm$ are hyperbolic buildings (with compact chambers).  For example, the building $I_{p,q+1}$ may be realised as a Kac--Moody building when $q$ is a power of a prime.  For further details, see Carbone--Garland~\cite{MR2017720} and R\'emy~\cite{MR1715140}. 

Apart from right-angled and Kac--Moody buildings, there are very few known constructions of hyperbolic buildings of dimension greater than $2$.  Haglund--Paulin~\cite{MR1957268} have constructed some three--dimensional hyperbolic buildings using ``tree-like" decompositions of the corresponding Coxeter systems, while Davis~\cite{MR2486800} gives covering-theoretic constructions of some three-dimensional hyperbolic buildings, including some where not all links are the same.  

As discussed in Section \ref{ss:coxeter} above, there is no complete classification of geometric hyperbolic Coxeter systems, and so there seems little hope of classifying general hyperbolic buildings.  Even for Fuchsian buildings, where the associated Coxeter systems are classified, the possible links $L$ are generalised $m$--gons, which have not been classified.  Indeed even after fixing a Coxeter system $(W,S)$ and the link $L$ at each vertex, there may be uncountably many hyperbolic buildings of type $(W,S)$ with links $L$ (see for example \cite[Theorem 3.6]{MR1877215}).  There are however cases in which ``local data" does determine the building.  For example, Bourdon's building $I_{p,q}$ is the unique simply-connected polygonal complex such that all faces are right-angled hyperbolic $p$--gons, and all links are $K_{q,q}$ \cite{MR1445387,MR1612350}.  For more on the question of uniqueness, see \cite[Section 2.3]{MR2807842}.

\section{Lattices}\label{s:lattices}

We now discuss lattices in hyperbolic buildings.  The state of the theory is such that for many hyperbolic buildings, the basic question of the existence of lattices in their automorphism groups is open.  We thus begin by describing known constructions of lattices in Section~\ref{ss:constructions}, then discuss a range of other questions concerning lattices in Section~\ref{ss:questions}.

\subsection{Constructions of lattices}\label{ss:constructions}

Let $\Delta$ be a hyperbolic building, and recall the characterisation of lattices in $\Aut(\Delta)$ from the introduction.  It will be seen that much more is known about constructions of cocompact lattices in $\Aut(\Delta)$ than about constructions of non-cocompact lattices.

If $\Delta$ is the universal cover of a finite polyhedral complex, as for example in the constructions of Fuchsian buildings due to Vdovina~\cite{MR1938699}, then the fundamental group of that finite polyhedral complex is a cocompact lattice in $\Aut(\Delta)$, since it acts freely and cocompactly on $\Delta$.     

Now suppose $\Delta$ is the universal cover of a complex of finite groups, over a finite underlying polyhedral complex $Y$.  For instance, one of the constructions of Bourdon's building $I_{p,q}$ is as the universal cover of a complex of finite groups over a right-angled hyperbolic $p$--gon $P$.  Let $\G$ be the fundamental group of this complex of groups.  Roughly speaking, $\G$ is an amalgam of the finite groups associated to the cells of $Y$.  Then $\G$ is a cocompact lattice in $\Aut(\Delta)$, since $\G$ acts on $\Delta$ with finite stabilisers and compact quotient $Y$.  

In \cite{Futer:2010uq}, Futer--Thomas constructed cocompact lattices in $\Aut(I_{p,q})$ which are fundamental groups of complexes of groups over $Y$ a (tessellated) surface.  They also showed that for some $p$ and $g$, whether there exists a cocompact lattice $\G < \Aut(I_{p,q})$ so that the quotient by the action of $\G$ is a genus $g$ surface depends upon the value of $q$.  This is the only known case where the values of the parameters $p$ and $q$ affect the existence of lattices in $\Aut(I_{p,q})$.

Complexes of finite groups may also be used to construct non-cocompact lattices.  In this case, the underlying complex $Y$ is infinite, and the assigned finite groups must have orders growing fast enough that the series in Equation \eqref{e:lattice} above converges.  For example, Thomas~\cite{MR2253444} obtained many cocompact and non-cocompact lattices for right-angled buildings by constructing a functor from graphs of groups,  with tree lattices as their fundamental groups, to complexes of groups, with right-angled building lattices as their fundamental groups.
More elaborate complexes of groups were used by Thomas to construct both cocompact and non-cocompact lattices for certain Fuchsian buildings in~\cite{MR2292989} and Davis--Moussong complexes in~\cite{Thomas:2008kx}.  

When $\Delta$ is a Davis--Moussong $(k,L)$--complex, then there is an associated Coxeter group $W(k,L)$, which is \emph{not} the type of the building $\Delta$, and the group $W(k,L)$ may be regarded as a cocompact lattice in $\Aut(\Delta)$.  In~\cite[Example 1.5.2]{MR1756974} Bourdon gives a lattice construction which ``lifts" lattices for Euclidean buildings to cocompact and non-cocompact lattices for certain Fuchsian buildings.    
 
Now suppose that the hyperbolic building $\Delta$ is the building for a Kac--Moody group $\Lambda$.  Recall that $\Lambda$ acts diagonally on the product $\Delta_+ \times \Delta_-$, where $\Delta_\pm \cong \Delta$.  Carbone--Garland~\cite{MR2017720} and independently R\'emy~\cite{MR1715140} showed that the stabiliser in $\Lambda$ of any point in the negative building $\Delta_-$ is a non-cocompact lattice in the automorphism group of the positive building $\Delta_+$.  These stabilisers may also be considered as lattices in the \emph{complete Kac--Moody group} $\hat{\Lambda} = \hat{\Lambda}_+$, which is a totally disconnected locally compact group acting on $\Delta_+$, and is obtained by completing $\Lambda$ using one of several methods.  Some of the lattices in $\hat\Lambda$ constructed in Gramlich--Horn--M\"uhlherr~\cite{MR2788084} are for hyperbolic buildings.  We do not know of any other constructions of lattices in complete Kac--Moody groups whose associated buildings are hyperbolic.

\subsection{Questions about lattices}\label{ss:questions}

As mentioned in the introduction, many questions concerning lattices in hyperbolic buildings are motivated by comparison with known results concerning lattices in semisimple Lie groups.  This background and motivation for the questions below is treated much more thoroughly in the corresponding sections of~\cite{MR2807842}, to which we refer the reader.  In general, a lot more is known about cocompact than about non-cocompact lattices in hyperbolic buildings.  There are also cases where the behaviour of compact and non-cocompact lattices is dramatically different.  It does appear that there is greater ``rigidity" when the hyperbolic building has (compact) simplicial chambers, perhaps because the associated Coxeter system then has the property, like all irreducible affine and finite Coxeter systems, that all $m_{st}$ are finite.

\subsubsection{Classification}

Once it is known that the automorphism group of a hyperbolic building $\Delta$ admits lattices, an immediate next question is to classify the lattices in $\Aut(\Delta)$ up to conjugacy.  This has only been done in special cases.  For instance in \cite{MR2665778}, Kangaslampi--Vdovina classify the torsion-free groups which act simply transitively on the vertices of Fuchsian buildings with triangular chambers and links the smallest generalised $4$--gon, and in \cite{Carbone:2011fk}, Carbone--Kangaslampi--Vdovina classify all such groups with torsion.

\subsubsection{Commensurability and commensurators}

Lattices may also be classified up to commensurability.  Recall that two subgroups $\G_1, \G_2 < G$ are \emph{commensurable} if there exists $g \in G$ so that $g\G_1 g^{-1} \cap \G_2$ has finite index in both $g \G_1 g^{-1}$ and $\G_2$.  Haglund \cite[Theorem 1.1]{MR2240922} proved that for $p \geq 6$, all cocompact lattices in $\Aut(I_{p,q})$ are commensurable.  In contrast, there are uncountably many commensurability classes of non-cocompact lattices in $\Aut(\Delta)$ for $\Delta$ a regular right-angled building (Thomas \cite[Main Theorem 2(b)]{MR2643689}).  

The \emph{commensurator} of a lattice $\G < G$ is the subgroup consisting of elements $g \in G$ such that $g \G g^{-1}$ and $\G$ are commensurable.  Haglund \cite{MR2413337} and independently Kubena--Thomas \cite{Barnhill:2008ys} proved that for $\Delta$ a regular right-angled building, the commensurator of a canonical cocompact lattice is dense in $G = \Aut(\Delta)$.   The question of commensurators of non-cocompact lattices is wide open, even for $I_{p,q}$.  

\subsubsection{Covolumes}

A basic question is to determine the set 
\[ \{ \mu(\G \backslash G) \mid \G < G \mbox{ is a lattice}\} \]
of covolumes of lattices in $G$.  Aspects of this question have been considered by Thomas for certain right-angled buildings \cite{MR2643689}, Fuchsian buildings \cite{MR2292989} and Davis--Moussong complexes \cite{Thomas:2008kx}, but many cases remain open.  For instance, it would be interesting to determine whether the set of covolumes of cocompact lattices in the Fuchsian buildings with triangular chambers considered in \cite{MR2665778} has a positive lower bound.

\subsubsection{Property (T) and finiteness properties}

Ballmann--\Swiatkowski~\cite{MR1465598}, Dymara--Januszkiewicz \cite{MR1946553} and \.Zuk \cite{MR1408975} have shown that the automorphism groups of many hyperbolic buildings $\Delta$ with simplicial chambers have Kazhdan's Property (T).  On the other hand, Corollary 3 of \cite{MR1465598} implies that the automorphism groups of Fuchsian buildings with chambers $p$--gons, $p\geq 4$, do not have Property (T).  In higher dimensions, it follows from Niblo--Reeves~\cite[Theorem B]{MR1432323} that the automorphism groups of right-angled buildings do not have Property (T), and Haglund--Paulin \cite[Theorem 1.5]{MR1668359} showed that the automorphism groups of hyperbolic buildings with chambers ``even" polytopes do not have Property (T).

A group $G$ has Property (T) if and only if all of its lattices have Property (T), and it is a well-known result of Kazhdan~\cite{MR0209390} that lattices with Property (T) are finitely generated.  All cocompact lattices in a hyperbolic building are finitely generated since they are fundamental groups of finite complexes of finite groups.  Infinite generation of some non-cocompact lattices for certain hyperbolic buildings was established by Thomas \cite{Thomas:2008kx} and Thomas--Wortman~\cite{MR2782548}.  It is not known whether for these buildings, there are any non-cocompact lattices which are finitely generated.

Very little is known about higher finiteness properties for lattices in hyperbolic buildings, apart from a recent result of Gandini~\cite{Gandini:2011vn} which bounds the homological finiteness length of non-cocompact lattices in $\Aut(X)$ for $X$ a locally finite contractible polyhedral complex.  As a corollary, such lattices are not finitely presentable.  The examples of \cite{Thomas:2008kx} and \cite{MR2782548} show that Gandini's bound is not sharp. 

\subsubsection{Residual finiteness, linearity and simplicity}

These questions are of particular interest for cocompact lattices in hyperbolic buildings with compact chambers, since such lattices are (finitely generated) word-hyperbolic groups.  Recall that a group $\G$ is \emph{residually finite} if for all $1 \neq \gamma \in \G$, there exists a finite quotient $\G \to \overline{\G}$ such that $\overline\gamma$ is nontrivial. It is a theorem of Mal'cev~\cite{MR0003420} that every finitely generated linear group is residually finite.   A major open conjecture of Gromov states that all word-hyperbolic groups are residually finite, while it is unknown whether every word-hyperbolic group is linear (see the introduction to Kapovich--Wise \cite{MR1735163}).  Word-hyperbolic groups are never simple \cite{MR919829, MR1357360}.

An important result, due to Wise \cite{MR1923477}, is the residual finiteness of cocompact lattices which are fundamental groups of complexes of finite groups over hyperbolic $p$--gons, for $p$ large enough (depending upon the angles of the polygon).  Combining this with \cite[Theorem 1.1]{MR2240922}, Haglund proved that for $p \geq 6$, all cocompact lattices in $\Aut(I_{p,q})$ are linear and thus residually finite.  The case $I_{5,q}$ is open.  As noted in Kangaslampi--Vdovina \cite{MR2665778}, the residual finiteness and linearity of the cocompact lattices in Fuchsian buildings with triangular chambers is also open.

For non-cocompact lattices, nonlinear examples in $\Aut(I_{p,q})$ were obtained by R\'emy~\cite{MR2084981} using Kac--Moody theory.  A striking recent result of Caprace--R\'emy \cite{MR2485882} is the simplicity of many Kac--Moody groups $\Lambda$, which are non-cocompact lattices in the product of their twin buildings $\Delta_+ \times \Delta_-$.

\subsubsection{Rigidity}

Various rigidity questions for lattices in hyperbolic buildings overlap with questions about the structure of the building itself and its boundary, and so we discuss these issues together here.  This section is only intended to be a list of some recent work in this area.

A careful proof of the folklore result that the visual boundary of a right-angled hyperbolic building is a Menger sponge was provided by Dymara--Osajda in \cite{MR2365885}.  Bourdon \cite[Theorems 1.1 and 1.2]{MR1445387} determined the conformal dimension of the visual boundary of $I_{p,q}$, and related this to its Hausdorff dimension.  A version of Mostow rigidity was also established by Bourdon \cite[Theorem 1.3]{MR1445387} for cocompact lattices in $I_{p,q}$, using combinatorial Patterson--Sullivan measures.  Bourdon--Pajot \cite{MR1789183} established quasi-isometric rigidity for $I_{p,q}$, and Xie \cite{MR2170496} generalised this result to all Fuchsian buildings.     

Since hyperbolic buildings are $\CAT(-1)$ spaces, several rigidity results for divergence groups apply to lattices in hyperbolic buildings.  As noted in the introduction to Burger--Mozes \cite{MR1325797}, the notion of a divergence group comes from Patterson--Sullivan theory for Kleinian groups.  If $\Aut(\Delta)$ acts cocompactly on the hyperbolic building $\Delta$, as is the case for most known constructions, then any nonelementary lattice $\G < \Aut(\Delta)$ is a divergence group \cite[Corollary 6.5(2)]{MR1325797}.  Hersonsky--Paulin \cite{MR1476054} generalised Mostow rigidity to divergence groups acting by isometries on many $\CAT(-1)$ spaces, including some hyperbolic buildings, and Burger--Mozes \cite{MR1325797} established $\CAT(-1)$ super-rigidity results for divergence groups.  

More recently, Daskalopoulos--Mese--Vdovina \cite{MR2827014} studied harmonic maps from symmetric spaces into target spaces including hyperbolic buildings, and as an application proved super-rigidity results for the isometry groups of a class of complexes including hyperbolic buildings.  Super-rigidity results for Kac--Moody groups $\Lambda$ were obtained by Caprace--R\'emy in \cite{MR2485882}.

An important ingredient in some classical rigidity results is the Howe--Moore property for unitary representations (see for example \cite{MR776417}), which concerns decay of matrix coefficients.  This property was shown not to hold for $\Aut(I_{p,q})$ by Bader--Shalom \cite[pp. 447--449]{MR2207022}, using Mackey theory.  

Finally, the volume entropy of hyperbolic buildings considers the asymptotic growth of volumes of balls in the building (by analogy with volume entropy for Riemannian manifolds).  This topic has been investigated by Hersonsky--Paulin \cite{MR1476054}, Leuzinger \cite{MR2248263} and most thoroughly by Ledrappier--Lim \cite{MR2643890}, using the geodesic flow on apartments and measures on suitable boundaries.

\bibliographystyle{siam}
\bibliography{refs} 

\end{document}